\documentclass[11pt]{article}
\usepackage[utf8]{inputenc}
\usepackage{graphicx}
\usepackage{caption}
\usepackage{amsmath}
\usepackage{amssymb}
\usepackage{float}
\usepackage{hyperref}
\usepackage{amsthm}
\usepackage{xcolor}
\usepackage{comment}
\usepackage{amsfonts}
\usepackage{dsfont}
\usepackage{stmaryrd}
\usepackage{xcolor}
\usepackage{subcaption}

\usepackage{siunitx}
\usepackage{mathtools} 
\mathtoolsset{showonlyrefs} 

\usepackage{lineno} 

\DeclareMathOperator{\R}{\mathbb{R}}

\def\Ec{{\cal E}}

\def\Fc{{\cal F}}
\def\Hc{{\cal H}}

\def\Pc{{\cal P}}
\def\Sc{{\cal S}}

\def\Wc{{\cal W}}
\def\Zc{{\cal Z}}
\def \E{\mathbb{E}}
\def \F{\mathbb{F}}
\def \M{\mathbb{M}}
\def \P{\mathbb{P}}

\def \H{\mathbb{H}}

\def\1{{\bf 1}}

\def \N{\mathbb{N}}

\newcommand{\di}{\mathrm{d}}

\def \Sum{\displaystyle\sum}

\def \b1{\bf{1}}

\def\boX{{\boldsymbol X}}

\def\boW{{\boldsymbol W}}
\def\boZ{{\boldsymbol Z}}

\def\bolx{{\boldsymbol x}}

\def\bosig{{\boldsymbol \sigma}}

\def \mry{\mathrm{y}} 
 
\def \mrz{\mathrm{z}}

\def\blue[#1]{{\color{blue}#1}}

\def\beqs{\begin{eqnarray*}}
\def\enqs{\end{eqnarray*}}
\def\beq{\begin{eqnarray}}
\def\enq{\end{eqnarray}}

\usepackage[algo2e,ruled,vlined]{algorithm2e} 
 \usepackage{algorithm}
 \usepackage{algorithmicx}
\usepackage{algpseudocode}
\usepackage{ marvosym }
\usepackage{csquotes}
\usepackage{hyperref}

\def\argmin{\mathop{\rm argmin}}
\def\argmin_#1{\underset{#1}{\mathrm{argmin\, }}}

\newtheorem{Theorem}{Theorem}[section]

\newtheorem{Assumption}[Theorem]{Assumption}

\newtheorem{Remark}[Theorem]{Remark}

\numberwithin{equation}{section} 

\usepackage[backend=bibtex, style=alphabetic]{biblatex}

\def \trans{^{\scriptscriptstyle{\intercal}}}

\usepackage{geometry}
\begin{document}
\title{Rate of convergence for particle approximation of \\ PDEs in Wasserstein space
\thanks{This work  was supported by FiME (Finance for Energy Market Research Centre) and
the ``Finance et D\'eveloppement Durable - Approches Quantitatives'' EDF - CACIB Chair.}}

\author{Maximilien \textsc{Germain}
\footnote{EDF R\&D, LPSM, Université de Paris  \sf \href{mailto:Maximilien.Germain at edf.fr}{mgermain at lpsm.paris}} \and  Huyên \textsc{Pham}
\footnote{LPSM, Université de Paris, and  FiME, and  CREST ENSAE \sf \href{mailto:pham at lpsm.paris}{pham at lpsm.paris}} \and  Xavier \textsc{Warin}
\footnote{EDF R\&D, FiME \sf \href{mailto:xavier.warin at  edf.fr}{xavier.warin at edf.fr}} 
}
\maketitle              
\begin{abstract}
We prove a rate of convergence for the $N$-particle approximation of 
a second-order partial diffe\-rential equation in the space of probability measures, like the Master equation or Bellman equation of mean-field control problem under common noise. 
The rate is of order $1/N$ for the pathwise error on the solution $v$ and of order $1/\sqrt{N}$ for the $L^2$-error on its $L$-derivative $\partial_\mu v$.  
The proof relies on backward stochastic differential equations techniques.  
\end{abstract}

\section{Introduction}

Let us consider the second-order parabolic partial differential equation (PDE) on the Wasserstein space $\Pc_2(\R^d)$ of square-integrable probability measures on $\R^d$, in the form: 
\begin{equation}  \label{PDEwasgen} 
\left\{ 
\begin{array}{rcl} 
 \partial_t v  +  \Hc(t,\mu,v,\partial_\mu v ,\partial_x\partial_\mu v,\partial_\mu^2 v) & = &    0, \quad \quad  (t,\mu) \in   [0,T)\times\Pc_2(\R^d),  \\
 v(T,\mu)   & =  &    G(\mu), \quad \quad \mu \in \Pc_2(\R^d). 
\end{array}
\right.
\end{equation}
Here, $\partial_\mu v(t,\mu)$ is the $L$-derivative on $\Pc_2(\R^d)$ (see \cite{cardel19}) of $\mu$ $\mapsto$ $v(t,\mu)$, and it is a function from $\R^d$ into $\R^d$,   
$\partial_x\partial_\mu v(t,\mu)$ is the usual derivative on $\R^d$ of $x$ $\in$ $\R^d$ $\mapsto$ $\partial_\mu v(t,\mu)(x)$ $\in$ $\R^d$,   hence valued in $\R^{d\times d}$ the set of $d\times d$-matrices with real coefficients, and 
$\partial_\mu^2 v(t,\mu)$ is the $L$-derivative of $\mu$ $\mapsto$ $\partial_\mu v(t,\mu)(.)$, hence a function from $\R^d\times\R^d$ into $\R^{d\times d}$. 
The terminal condition is given by a real-valued function $G$ on $\Pc_2(\R^d)$,  and the Hamiltonian $\Hc$ of this PDE is assumed to be in semi-linear (non linear w.r.t. $v$, $\partial_\mu v$,  and linear w.r.t. $\partial_x\partial_\mu$, $\partial_\mu^2 v$) expectation form: 
\begin{align} \label{Hintegral} 
\Hc(t,\mu,y,\mrz(.),\gamma(.),\gamma_0(.,.)) &=   \int_{\R^d} \big[ H(t,x,\mu,y,\mrz(x)) + \frac{1}{2}{\rm tr}\big((\sigma\sigma\trans+\sigma_0\sigma_0\trans)(t,x,\mu)\gamma(x)\big) \big] \mu(\di x) \\
& \quad  \quad  +  \;  \frac{1}{2} \int_{\R^d\times\R^d} {\rm tr}\big(\sigma_0(t,x,\mu)\sigma_0\trans(t,x',\mu)\gamma_0(x,x')\big) \mu(\di x)\mu(\di x'),  
\end{align}
for some real-valued measurable  function $H$ defined   on $[0,T]\times\R^d\times\Pc_2(\R^d)\times\R\times\R^d$, and where $\sigma$, $\sigma_0$ are measurable functions on $[0,T]\times\R^d\times\Pc_2(\R^d)$, valued respectively in
$\R^{d\times n}$, and $\R^{d\times m}$. 
Here ${\rm tr}(M)$ denotes the trace of a square matrix $M$, while $M\trans$ is its  transpose, and $.$ is the scalar product.   


\vspace{1mm}

PDEs in Wasserstein space have been largely studied in the literature over the last years, notably with the emergence of the mean-field game theory, and we mention among others the papers \cite{benetal15},  \cite{ganswi15}, \cite{phawei18}, \cite{caretal19}, \cite{saporitozhang19}, \cite{buretal20}, and other 
references in the two-volume monographs \cite{cardel19}-\cite{cardel2}.   

An important  application concerns mean-field type control  problems with common noise. The controlled stochastic McKean-Vlasov dynamics  is given by 
\begin{align} \label{dynMKVcontrol}
\di X_s^\alpha &= \; \beta(s,X_s^\alpha,\P^0_{_{X_s^\alpha}},\alpha_s) \di s + \sigma(s,X_s^\alpha,\P^0_{_{X_s^\alpha}}) \di W_s \\
& \quad \quad \quad + \;  \sigma_0(s,X_s^\alpha,\P^0_{_{X_s^\alpha}}) \di W^0_s ,  \quad t \leq s \leq T, \; X_t^\alpha =  \xi,  
\end{align}
 where  $W$ is  a $n$-dimensional Brownian motion, independent of a $m$-dimensional Brownian motion $W^0$ (representing the common noise) on a filtered probability space $(\Omega,\Fc,\F=(\Fc_t)_{0\leq t\leq T},\P)$, 
 the control process $\alpha$ is $\F$-adapted valued in some Polish space $A$, and here $\P^0$ denotes the conditional law given $W^0$. 
The value function defined on $[0,T]\times\Pc_2(\R^d)$  by 
\beqs
v(t,\mu) & = & \inf_{\alpha} \E_{t,\mu} \Big[ \int_t^T e^{-r(s-t)} f(X_s^\alpha,\P^0_{_{X_s^\alpha}},\alpha_s) \di s  + e^{-r(T-t)} g(X_T^\alpha,\P^0_{_{X_T^\alpha}}) \Big], 
\enqs
(here $\E_{t,\mu}[\cdot]$ is the conditional expectation  given that the law at time $t$ of $X$ solution to \eqref{dynMKVcontrol}  is equal  to  $\mu$) 
is shown to satisfy the Bellman  equation \eqref{PDEwasgen}-\eqref{Hintegral}  (see \cite{benetal13}, \cite{cospha19}, \cite{djepostan19})  with   $G(\mu)$ $=$  $\int g(x,\mu) \mu(dx)$, 
$\sigma$, $\sigma_0$ as in \eqref{dynMKVcontrol} and 
\begin{align}  \label{defhcontrol}
H(t,x,\mu,y,z) &=\;  -ry + \inf_{a \in A} \big[ \beta(t,x,\mu,a).z   + f(x,\mu,a) \big].  
\end{align}

\vspace{1mm}

We now consider a finite-dimensional approximation of the PDE  \eqref{PDEwasgen}-\eqref{Hintegral} in the Wasserstein space. This can be derived formally by  looking at the PDE for $\mu$ to ave\-rages of Dirac masses, and it turns out that the corresponding PDE takes the  form 
\begin{equation} \label{PDEfini} 
\begin{cases}
& \partial_t v^N + \frac{1}{N} \Sum_{i=1}^N H\big(t,x_i, \bar\mu(\bolx),v^N,N D_{x_i} v^N) + \frac{1}{2} {\rm tr}(\Sigma_N(t,\bolx) D_{\bolx}^2 v^N) \; = \;  0,   \mbox{ on }   [0,T) \times  (\R^d)^N \\
& v^N(T,\bolx) \; = \; G\big(\bar\mu(\bolx)\big), \quad  \bolx = (x_i)_{i \in \llbracket1,N \rrbracket}   \in  (\R^d)^N,   
\end{cases}
\end{equation}
where $\bar\mu(.)$ is the empirical measure function defined by $\bar\mu(\bolx)$ $=$ $\frac{1}{N}\sum_{i=1}^N \delta_{x_i}$, for any $\bolx$ $=$ $(x_1,\ldots,$ $x_N)$, $N$ $\in$ $\N^*$, and $\Sigma_N$ $=$ $(\Sigma_N^{ij})_{i,j\in \llbracket 1,N\rrbracket}$  is the $\R^{Nd\times Nd}$-valued  function with 
block matrices  $\Sigma_N^{ij}(t,\bolx)$ $=$ $\sigma(t,x_i,\bar\mu(\bolx))\sigma\trans(t,x_j,\bar\mu(\bolx))\delta_{ij}$ $+$ $\sigma_0(t,x_i,\bar\mu(\bolx))\sigma_0\trans(t,x_j,\bar\mu(\bolx))$ $\in$ $\R^{d\times d}$.  
In the special case where $H$ has the form \eqref{defhcontrol}, we notice that 
\eqref{PDEfini} is the Bellman equation for the $N$-cooperative problem,  whose convergence to the mean-field control problem  has been studied in \cite{lac17}, \cite{cardel2}, \cite{laurire2019backward,laurire2020convergence}, when  $\sigma_0$ $\equiv$ $0$  (no common noise), and 
recently by \cite{dje20} in the common noise case. We point out that these works do not consider the same master equation. In particular their master equation is stated on $[0,T]\times\R^d\times\Pc_2(\R^d) $ and is linear in $\partial_\mu u$ whereas we allow a non-linear dependence in this derivative. Moreover our master equation is in expectation form. In \cite{laurire2020convergence} the master equation is approached by a system of $N$ coupled PDEs on $[0,T]\times (\R^d)^N$ whereas we consider a single approximating PDE on $[0,T]\times (\R^d)^N$.
For more general Hamiltonian functions $H$,  it has been recently proved in \cite{ganmayswi20} that the sequence of viscosity solutions $(v^N)_N$ to  
\eqref{PDEfini} converge locally uniformly to the viscosity solution $v$  to \eqref{PDEwasgen} when $\sigma=0$ and $\sigma_0$ does not depend on space and measure arguments. For a detailed comparison between this work and ours, we refer to Remark \ref{rem: comparison gangbo mayorga swiech}.

\vspace{1mm}

In this paper, we adopt a probabilistic approach by considering a backward stochastic differential equation (BSDE) representation for the finite-dimensional PDE \eqref{PDEfini}  according to the classical work \cite{pardoux1990adapted}. 
The solution $(Y^N, \boZ^N=(Z^{i,N})_{1\leq i\leq N})$  to this BSDE is written with an underlying forward particle system $\boX^N=(X^{i,N})_{1\leq i\leq N}$ of a McKean-Vlasov SDE, and connected to the PDE  \eqref{PDEfini}  via the Feynman-Kac formula: 
$Y_t^N$ $=$ $v^N(t,\boX_t^N)$, $Z_t^{i,N}$ $=$ $D_{x_i} v^N(t,\boX_t^{i,N})$, $0\leq t\leq T$.  
By using BSDE techniques, our main contribution is  to show a rate of convergence of order $1/N$ of $|Y^N_t - v(t,\bar\mu(\boX_t^N))|$, 
and also of $|NZ_t^{i,N} - \partial_\mu v(t,\bar\mu(\boX_t^N))(X_t^{i,N})|^2$, 
$i$ $=$ $1,\ldots,N$, for suitable norms, and under some regularity conditions on  $v$ (see Theorem \ref{th: N convergence} and Theorem \ref{th: convergence Z}). 
This rate of convergence on the particles approximation of  $v$ and its $L$-derivative  is new to the best of  our knowledge. 
We point out that classical BSDE arguments for proving the rate of convergence do not apply directly due to the presence of the factor $N$ in front of $D_{x_i} v^N$ in the generator $H$, and we rather use linearization arguments and change of probability measures to overcome these issues. Another issue is due to the fact that the BSDE dimension $d\times N$ is exploding with the number of particles therefore we have to track down the influence of the dimension in the estimations, whereas classical BSDE works usually consider a fixed dimension $d$ which is incorporated into constants. 
 
\vspace{1mm}

The outline of the paper is organized as follows. In Section \ref{secparticle}, we formulate the particle approximation of the PDE and its BSDE representation, and state the rate of convergence for $v$ and its $L$-derivative. 
Section \ref{secproof} is devoted to the proof of these results.

\section{Particles approximation of  Wasserstein PDEs} \label{secparticle}

The formal derivation of the finite-dimensional approximation PDE is obtained as follows. We look at the PDE  \eqref{PDEwasgen}-\eqref{Hintegral}  for $\mu$ $=$ $\bar\mu(\bolx)$ $=$ $\frac{1}{N}\sum_{i=1}^N \delta_{x_i}$ $\in$ $\Pc_2(\R^d)$, 
when $\bolx$ $=$  $(x_i)_{i\in \llbracket 1,N\rrbracket}$ runs over  $(\R^d)^N$. By setting  $\tilde v^N(t,\bolx)$ $=$ $v(t,\bar\mu(\bolx))$, and assuming that $v$ is smooth, we  have for all $(i,j)$ $\in$ $\llbracket 1, N \rrbracket$  (see Proposition 5.35 and Proposition 5.91 in \cite{cardel19}): 
 \begin{equation} \label{derivN} 
 \begin{cases}
& D_{x_i} \tilde v^N(t,\bolx)  \; = \;  \frac{1}{N} \partial_\mu v(t,\bar\mu(\bolx))(x_i),  \\
& D^2_{x_i x_j} \tilde v^N(t,\bolx)  \; = \;  \frac{1}{N} \partial_x \partial_\mu v(t,\bar\mu(\bolx))(x_i)\delta_{ij} + \frac{1}{N^2} \partial^2_\mu  v(t,\bar\mu(\bolx))(x_i,x_j). 
\end{cases}
\end{equation}
 By  substituting into the PDE  \eqref{PDEwasgen}-\eqref{Hintegral}  for $\mu$ $=$ $\bar\mu(\bolx)$, and using \eqref{derivN}, we then see that $\tilde v^N$ satisfies the relation: 
\begin{align} \label{tildevN} 
\partial_t \tilde v^N +  \frac{1}{N} \Sum_{i=1}^N H\big(t,x_i, \bar\mu(\bolx),\tilde v^N,N D_{x_i} \tilde v^N) & \\
+ \; \frac{1}{2} \sum_{i=1}^N {\rm tr}\big[ \big(\sigma\sigma\trans + \sigma_0\sigma_0\trans\big)(t,x_i,\bar\mu(\bolx)) \big(D_{x_i}^2 \tilde v^N - \frac{1}{N^2} \partial_\mu^2 v(t,\bar\mu(\bolx)) (x_i,x_i) \big) \big]  & \\
+ \;  \frac{1}{2} \sum_{i\neq j \in \llbracket 1,N\rrbracket} {\rm tr} \big( \sigma_0(t,x_i,\bar\mu(\bolx))\sigma_0\trans(t,x_j,\bar\mu(\bolx)) D_{x_ix_j}^2 \tilde v^N \big) & \\
+ \; \frac{1}{2N^2} \sum_{i=1}^N {\rm tr}\big(\sigma_0\sigma_0\trans(t,x_i,\bar\mu(\bolx))  \partial_\mu^2 v(t,\bar\mu(\bolx)) (x_i,x_i) \big) & = \; 0
\end{align} 
for $(t,\bolx=(x_i)_{i\in\llbracket 1,N\rrbracket})$ $\in$ $[0,T)\times(\R^d)^N$, together with the terminal condition $\tilde v^N(t,\bolx)$ $=$ $G(\bar\mu(\bolx))$.  
By neglecting the terms  $\partial_\mu^2 v/N^2$ in the above relation, we obtain the  PDE \eqref{PDEfini} for $v^N$ $\simeq$ $\tilde v^N$. 
The purpose of this section is to  rigorously justify this approximation and state a rate of convergence for $v^N$ towards $v$, as well as a convergence for their gradients.

 \subsection{Particles BSDE approximation}

Let us  introduce an arbitrary measurable $\R^d$-valued function $b$ on $[0,T]\times\R^d\times\Pc_2(\R^d)$, and set $B_N$ the $(\R^d)^N$-valued function defined on $[0,T]\times(\R^d)^N$ by $B_N(t,\bolx)$ $=$ $(b(t,x_i,\bar\mu(\bolx))_{i\in \llbracket 1,N\rrbracket}$ for $(t,\bolx=(x_i)_{i\in\llbracket 1,N\rrbracket})$ $\in$ 
$[0,T)\times(\R^d)^N$. The finite-dimensional PDE  \eqref{PDEfini} may then be written as 
\begin{equation} \label{PDEfinib} 
\begin{cases}
& \partial_t v^N + B_N(t,\bolx).D_{\bolx} v^N + \frac{1}{2}{\rm tr}\big(\Sigma_N(t,\bolx)D_{\bolx}^2 v^N\big) \\
& \quad \quad + \;   \frac{1}{N} \Sum_{i=1}^N H_b\big(t,x_i, \bar\mu(\bolx),v^N,N D_{x_i} v^N)    \; = \;  0, \quad \mbox{ on }   [0,T) \times  (\R^d)^N,  \\
& v^N(T,\bolx) \; = \; G\big(\bar\mu(\bolx)\big), \quad  \bolx = (x_i)_{i \in \llbracket1,N \rrbracket}   \in  (\R^d)^N,   
\end{cases}
\end{equation}
where $H_b(t,x,\mu,y,z)$ $:=$ $H(t,x,\mu,y,z)$ $-$ $b(t,x,\mu).z$.  For error analysis purpose, the function $b$ can be simply taken to be zero. The introduction of the function $b$ is actually motivated by numerical purpose. It corresponds indeed to the drift of training simulations for approximating the function $v^N$, 
notably by machine learning methods, and should be chosen for suitable exploration of the state space, see  a detailed discussion in  our companion paper \cite{gerlauphawar21a}. In this paper, we fix an arbitrary function $b$ (satisfying Lipschitz condition to be precised later).

Following \cite{pardoux1990adapted}, it is well-known that the semi-linear PDE \eqref{PDEfinib}  admits a probabilistic representation in terms  of forward backward SDE. The forward component  is defined by the process 
$\boX^N$ $=$ $(X^{i,N})_{i\in\llbracket 1,N\rrbracket}$ valued in $(\R^d)^N$, solution to the SDE: 
\begin{align} \label{diffboX}
\di\boX_t^N &= \; B_N(t,\boX_t^N) \di t    + \sigma_N(t,\boX_t^N) \di \boW_t + \bosig_0(t,\boX_t^N) \di W^0_t
\end{align}
where $\sigma_N$ is the block diagonal matrix with block diagonals $\sigma_N^{ii}(t,\bolx)$ $=$ $\sigma(t,x_i,\bar\mu(\bolx))$, 
$\bosig_0$ $=$ $(\sigma_0^i)_{i\in\llbracket 1,N\rrbracket}$ is the $(\R^{d\times m})^N$-valued function with 
$\bosig_0^{i}(t,\bolx)$ $=$ $\sigma_0(t,x_i,\bar\mu(\bolx))$, 
for $\bolx$ $=$ $(x_i)_{i\in\llbracket 1,N\rrbracket}$, $\boW$ $=$ $(W^1,\ldots,W^N)$ where $W^i$, $i$ $=$ $1,\ldots,N$, are independent $n$-dimensional Brownian motions, independent of a $m$-dimensional Brownian motion $W^0$ on a filtered probability space $(\Omega,\Fc,\F=(\Fc_t)_{0\leq t\leq T},\P)$. 
Notice that $\Sigma_N$ $=$ $\sigma_N\sigma_N\trans$ $+$ $\bosig_0\bosig_0\trans$, and  $\boX^N$ is the particles system of the McKean-Vlasov SDE: 
\begin{align} \label{SDEMKV} 
\di X_t &= \; b(t,X_t,\P_{_{X_t}}) \di t + \sigma(t,X_t,\P^0_{_{X_t}}) \di W_t + \sigma_0(t,X_t,\P^0_{_{X_t}}) \di W_t^0,  
\end{align}
where $W$ is an $n$-dimensional Brownian motion independent of $W^0$. 
The backward component is defined by the pair process $(Y^N,\boZ^N=(Z^{i,N})_{i\in\llbracket 1,N\rrbracket})$ valued in $\R\times(\R^d)^N$, solution to
\begin{align} \label{BSDEMKVN}
Y_t^N & = \; G\big(\bar\mu(\boX_{T}^N)\big) + \frac{1}{N} \sum_{i=1}^N \int_t^T H_b(s,X_s^{i,N},\bar\mu(\boX_{s}^N) ,Y_s^N, N Z_s^{i,N}) \di s\\
& \quad \quad  - \; \sum_{i=1}^N \int_t^T (Z_s^{i,N})\trans \sigma\big(s,X_s^{i,N},\bar\mu(\boX_{s}^N) \big)  \di W_s^i, \\
& \quad \quad  - \; \sum_{i=1}^N \int_t^T (Z_s^{i,N})\trans \sigma_0\big(s,X_s^{i,N},\bar\mu(\boX_{s}^N) \big)  \di W_s^0, \quad \quad   0 \leq t \leq T. 
\end{align}
We shall assume that the measurable functions  $(t,x,\mu)$ $\mapsto$ $b(t,x,\mu)$, $\sigma(t,x,\mu)$ satisfy a Lipschitz condition in $(x,\mu)$ $\in$ $\R^d\times\Pc_2(\R^d)$ uniformly w.r.t. $t$ $\in$ $[0,T]$, which ensures the existence and uniqueness of a strong solution $\boX^N$ 
$\in$ $\Sc_\F^2((\R^d)^N)$ 
to \eqref{diffboX} given an initial condition. Here, $\Sc_\F^2(\R^q)$ is the set of $\F$-adapted process $(V_t)_t$ valued in $\R^q$ s.t. $\E\big[\sup_{0\leq t\leq T}|V_t|^2\big]$ $<$ $\infty$,  ($|.|$ is the Euclidian norm on $\R^q$, and for a matrix $M$, we choose the Frobenius norm $|M|$ $=$ $\sqrt{{\rm tr}(MM\trans)}$) and  the Wasserstein space $\Pc_2(\R^d)$ is endowed with the Wasserstein distance 
\begin{align}
\Wc_2(\mu,\mu') &= \;  \Big( \inf\big\{ \E|\xi -\xi'|^2: \xi \sim \mu, \xi' \sim \mu' \big\} \Big)^{\frac{1}{2}},
\end{align}
and we set $\|\mu\|_{_2}$ $:=$ $\big( \int_{\R^d} |x|^2\ \mu(\di x)\big)^{\frac{1}{2}}$ for $\mu$ $\in$ $\Pc_2(\R^d)$. Assuming also that the measurable function $(t,x,\mu,y,z)$ $\mapsto$ $H_b(t,x,\mu,y,z)$  is Lipschitz in $(y,z)$ $\in$ $\R\times\R^d$ uniformly with respect to 
$(t,x,\mu)$ $\in$ $[0,T]\times\R^d\times\Pc_2(\R^d)$, and the measurable function $G$ satisfies a quadratic growth condition on $\Pc_2(\R^d)$, we have  the existence and uniqueness of a solution $(Y^N,\boZ^N=(Z^{i,N})_{i\in\llbracket 1,N\rrbracket})$   $\in$ $\Sc_\F^2(\R)\times\H_\F^2((\R^d)^N)$   to \eqref{BSDEMKVN}, and the  connection with the PDE \eqref{PDEfinib}  (satisfied in general in the viscosity sense) via the (non linear) Feynman-Kac formula: 
\begin{align} \label{relYNv}
Y_t^N &= \; v^N(t,\boX_t^N),  \quad \mbox{ and } \quad  Z_t^{i,N} \; = \;  D_{x_i} v^N(t,\boX_t^N), \;\; i=1,\ldots,N, \;\;  0 \leq t \leq T,
\end{align}
(when $v^N$ is smooth for the last relation).  Here, $\H_\F^2(\R^q)$ is the set of $\F$-adapted process $(V_t)_t$ valued in $\R^q$ s.t.  $\E\big[ \int_0^T |V_t|^2 \di t\big]$ $<$ $\infty$. 

\subsection{Main results}

We aim to analyze the particles approximation error on the solution $v$ to the PDE \eqref{PDEwasgen}, and its $L$-derivative $\partial_\mu v$  
by considering the pathwise error on $v$:
\begin{align}
\Ec_N^\mry & := \;  \sup_{0\leq t \leq T} \big|Y_t^N - v(t,\bar\mu(\boX_t^N)) \big|,  
\end{align}
and the $L^2$-error on its $L$-derivative
\begin{align}
\big\| \Ec_N^\mrz \big\|_{_2}  & := \; \frac{1}{N} \sum_{i=1}^N \Big( \int_0^T \E\big|NZ_t^{i,N} -   \partial_\mu  v(t,\bar\mu(\boX_t^N))(X_t^{i,N}) \big|^2 \di t \Big)^{\frac{1}{2}}, 
\end{align}
where the initial conditions of the particles system, $X_0^{i,N}$, $i$ $=$ $1,\ldots,N$, are i.i.d. with distribution $\mu_0$.

 \vspace{1mm}

Here, it is assumed that we have the existence and uniqueness of a classical solution $v$ to the PDE \eqref{PDEwasgen}-\eqref{Hintegral}. More precisely, we make the following assumption:

\begin{Assumption}[Smooth solution to the Master Bellman PDE]\label{assumption: solution}
There exists a unique solution $v$ to  \eqref{PDEwasgen}, which lies in  $C_b^{1,2}([0,T]\times\Pc_2(\R^d))$  that is:
\begin{itemize}
\item $v(.,\mu)$ $\in$ $C^1([0,T))$, and continuous on $[0,T]$, for any $\mu$ $\in$ $\Pc_2(\R^d)$,  
\item $v(t,.)$ is fully $C^2$ on $\Pc_2(\R^d)$, for any $t$ $\in$ $[0,T]$  in the sense that: $(x,\mu)$ $\in$ $\R^d\times\Pc_2(\R^d)$ $\mapsto$ $\partial_\mu v(t,\mu)(x)$ $\in$ $\R^d$,  $(x,\mu)$ $\in$ $\R^d\times\Pc_2(\R^d)$ $\mapsto$ $\partial_x \partial_\mu v(t,\mu)(x)$ $\in$ $\M^d$, and 
$(x,x',\mu)$ $\in$ $\R^d\times\R^d$ $\times\Pc_2(\R^d)$ $\mapsto$ $\partial_\mu^2 v(t,\mu)(x,x')$ $\in$ $\M^d$, are well-defined and jointly continuous, 
\item there exists some constant $L$ $>$ $0$ s.t. for all $(t,x,\mu)$ $\in$ $[0,T]\times\R^d\times\Pc_2(\R^d)$
\begin{align}
\big| \partial_\mu v(t,\mu)(x) \big| & \leq \; L\big(1 + |x| + \|\mu\|_{_2} \big), \quad \big| \partial_\mu^2 v(t,\mu)(x,x) \big| \; \leq \; L. 
\end{align}
\end{itemize}
\end{Assumption}

The existence of classical solutions to mean-field PDE in Wasserstein space is a cha\-llenging problem, and beyond the scope of this paper. We refer to \cite{buckdahn2017}, \cite{chacridel15}, \cite{saporitozhang19}, \cite{wuzhang20}  
for conditions ensuring regularity results of some Master PDEs. Notice also that linear-quadratic mean-field control problems have explicit smooth solutions as in Assumption \ref{assumption: solution}, see e.g. \cite{phawei18}.

\vspace{1mm}

We also make some rather standard assumptions on the coefficients of the forward backward SDE:

\begin{Assumption}[Lipschitz condition on the coefficients of the forward backward SDE]\label{assumption: Lipschitz}
\noindent 
\begin{itemize}
\item[(i)] The drift and volatility coefficients $b,\sigma,\sigma_0$ are Lipschitz: there exist positive constants $[b]$, $[\sigma]$, and $[\sigma_0]$ s.t. for all $t$ $\in$ $[0,T]$, $x,x'$ $\in$ $\R^d$, $\mu,\mu'$ $\in$ $\Pc_2(\R^d)$, 
\begin{align}
 |b(t,x,\mu) - b(t,x',\mu')| & \leq \;  [b] \big(|x-x'| + \Wc_2(\mu,\mu')\big) \\
 |\sigma(t,x,\mu) - \sigma(t,x',\mu')| &\leq \;  [\sigma] \big(|x-x'| + \Wc_2(\mu,\mu')\big) \\
 |\sigma_0(t,x,\mu) - \sigma_0(t,x',\mu')| &\leq \;  [\sigma_0] \big(|x-x'| + \Wc_2(\mu,\mu')\big).
\end{align}
\item[(ii)] For all $(t,x,\mu)$ $\in$ $[0,T]\times\R^d\times\Pc_2(\R^d)$, $\Sigma(t,x,\mu)$ $:=$ $\sigma\sigma\trans(t,x,\mu)$ is invertible, and the function $\sigma$, and its pseudo-inverse $\sigma^{+}$ $:=$ $\sigma\trans\Sigma^{-1}$ are bounded. 
\item[(iii)] $\mu_0$ $\in$ $\Pc_{4q}(\R^d)$ for some $q$ $>$ $1$, i.e., $\|\mu_0\|_{_{4q}}$ $:=$ $\big(\int |x|^{4q}\mu_0(\di x))^{\frac{1}{4q}}$ $<$ $\infty$, and 
\begin{align}
 \int_0^T |b(t,0,\delta_0)|^{4q} +   |\sigma(t,0,\delta_0)|^{4q}  +  |\sigma_0(t,0,\delta_0)|^{4q} \; \di t & < \; \infty. 
 \end{align}
\item[(iv)] The driver $H_b$ satisfies the  Lipschitz condition: there exist positive constants $[H_b]_{_1}$ and $[H_b]_{_2}$ s.t. for all $t$ $\in$ $[0,T]$, $x,x'$ $\in$ $\R^d$, $\mu,\mu'$ $\in$ $\Pc_2(\R^d)$, $y,y'$ $\in$ $\R$, $z,z'$ $\in$ $\R^d$, 
\begin{align}
| H_b(t,x,\mu,y,z) - H_b(t,x,\mu,y',z')|  & \leq \;   [H_b]_{_1} ( |y-y'| + |z-z'|) \label{lipHb} \\
| H_b(t,x,\mu,y,z) - H_b(t,x',\mu',y,z)| &\leq \;   [H_b]_{_2} \big(1+|x|+|x'|+  \|\mu\|_{_2} + \|\mu'\|_{_2}\big)  \nonumber \\
&  \quad \quad   \big( |x-x'| + \Wc_2(\mu, \mu')\big). \label{lipHb2}
\end{align}
\item[(v)] The terminal condition satisfies the (locally)  Lipschitz condition: there exists some positive constant $[G]$ s.t. for all $\mu,\mu'$ $\in$ $\Pc_2(\R^d)$
\begin{align}
|G(\mu) - G(\mu')| &\leq \;  [G] \big(\|\mu\|_{_2} + \|\mu'\|_{_2} \big) \Wc_2(\mu,\mu').
\end{align} 
\end{itemize}
\end{Assumption}

\vspace{1mm}

In order to have a convergence result for the first order Lions derivative we have to make a stronger assumption.

\begin{Assumption}\label{assumption: linear}
\noindent 
\begin{itemize}
    \item[(i)]  The function $H_b$ is in the form: 
    \begin{equation}
    H_b(t,x,\mu,y,z) = H_1(t,x,\mu,y) + H_2(t,\mu,y).z,
    \end{equation} 
where  $H_1: [0,T]\times \R^d\times\Pc_2(\R^d)\times \R \mapsto \R$ verifies for all $t$ $\in$ $[0,T]$, $x,x'$ $\in$ $\R^d$, $\mu,\mu'$ $\in$ $\Pc_2(\R^d)$, $y,y'$ $\in$ $\R$, $z,z'$ $\in$ $\R^d$, 
\begin{align}
| H_1(t,x,\mu,y) - H_1(t,x,\mu,y')|  & \leq \;   [H_1]_{_1} |y-y'|  \label{lipH1} \\
| H_1(t,x,\mu,y) - H_1(t,x',\mu',y)| &\leq \;   [H_1]_{_2} \big(1+|x|+|x'|+  \|\mu\|_{_2} + \|\mu'\|_{_2}\big)  \nonumber \\
&  \quad \quad   \big( |x-x'| + \Wc_2(\mu, \mu')\big), 
\end{align} and $H_2: [0,T]\times \R^d\times\Pc_2(\R^d)\times \R \mapsto \R^d$ is bounded and verifies for all $t$ $\in$ $[0,T]$, $x,x'$ $\in$ $\R^d$, $\mu,\mu'$ $\in$ $\Pc_2(\R^d)$, $y,y'$ $\in$ $\R$, $z,z'$ $\in$ $\R^d$, 
\begin{align}
| H_2(t,x,\mu,y) - H_2(t,x,\mu,y')|  & \leq \;   [H_2]_{_1} |y-y'|  \label{lipH2} \\
| H_2(t,x,\mu,y) - H_2(t,x',\mu',y)| &\leq \;   [H_2]_{_2} \big(1+|x|+|x'|+  \|\mu\|_{_2} + \|\mu'\|_{_2}\big)  \nonumber \\
&  \quad \quad   \big( |x-x'| + \Wc_2(\mu, \mu')\big). 
\end{align}
    \item[(ii)]  $\sigma_0$ is uniformly elliptic and does not depend on $x$, namely there  exists $c_0>0$ such that for all $t\in[0,T]$, $\mu\in\Pc_2(\R^d)$, $z\in\R^d$
    \begin{equation}
       z\trans \sigma_0(t,\mu) \sigma_0\trans(t,\mu) z \; \geq \;  c_0  |z|^2.
    \end{equation}
    \item[(iii)]  There exists some constant $L$ $>$ $0$ s.t. for all $(t,x,\mu)$ $\in$ $[0,T]\times\R^d\times\Pc_2(\R^d)$
\begin{align}
\big| \partial_\mu v(t,\mu)(x) \big| & \leq \; L. 
\end{align}
\end{itemize}
\end{Assumption} 

\begin{Remark} \label{remlip}
The Lipschitz condition on $b$, $\sigma$ in Assumption \ref{assumption: Lipschitz}(i) implies  that the functions $\bolx$ $\in$ 
$(\R^d)^N$  $\mapsto$ $B_N(t,\bolx)$, resp. $\sigma_N(t,\bolx)$ and $\bosig_0(t,\bolx)$, defined in \eqref{diffboX}, 
are Lipschitz (with Lipschitz constant $2[b]$, resp. $2[\sigma]$ and  $2[\sigma_0]$). Indeed, we have 
\begin{align}
|B_N(t,\bolx) - B_N(t,\bolx')|^2 & = \sum_{i=1}^N |b(t,x_i,\bar\mu(\bolx)) - b(t,x_i,\bar\mu(\bolx'))|^2\\
    & \leq 2[b]^2 \sum_{i=1}^N (|x_i - x_i'|^2 + \Wc_2(\bar\mu(\bolx),\bar\mu(\bolx'))^2)\\
    & \leq 2[b]^2  (|\bolx-\bolx'|^2 + \sum_{i=1}^N  \frac{1}{N}|\bolx-\bolx'|^2) \; = \; 4[b]^2|\bolx-\bolx'|^2,
\end{align}
for $\bolx$ $=$ $(x_i)_{i\in\llbracket 1,N\rrbracket}$, and similarly for $\sigma_N$ and $\bosig_0$. 

This yields the existence and uniqueness of a solution 
$\boX^N$ $=$ $(X^{i,N})_{i\in\llbracket 1,N\rrbracket}$ to  \eqref{diffboX} given initial conditions. Moreover, under Assumption \ref{assumption: Lipschitz}(iii), we have the standard estimate: 
\begin{align} \label{estiXiN} 
 \E\big[ \sup_{0\leq t\leq T} |\boX_t^{N}|^{4q} \big] & \leq \; C \big( 1 + \|\mu_0\|^{4q}_{_{4q}} \big) \; < \; \infty, \quad i=1,\ldots,N, 
 \end{align}
for some constant $C$ (possibly depending on $N$). The Lipschitz condition on $H_b$ w.r.t. $(y,z)$ in Assumption \ref{assumption: Lipschitz}(iv), and the quadratic growth condition on $G$ from  Assumption \ref{assumption: Lipschitz}(v) gives the existence and uniqueness of a solution 
$(Y^N,\boZ^N=(Z^{i,N})_{i\in\llbracket 1,N\rrbracket})$   $\in$ $\Sc_\F^2(\R)\times\H_\F^2((\R^d)^N)$   to \eqref{BSDEMKVN}. Moreover, by Assumption \ref{assumption: Lipschitz}(iv)(v), we see that 
\begin{align}
 &  \big| \frac{1}{N} \sum_{i=1}^N H_b(t,x_i,\bar\mu(\bolx),y,z_i) - \frac{1}{N} \sum_{i=1}^N H_b(t,x_i',\bar\mu(\bolx'),y,z_i)\big| \\
\leq & \;    [H_b]_{_2} \frac{1}{N} \sum_{i=1}^N  \big(1+|x_i|+|x_i'|  + \frac{1}{\sqrt{N}}(|\bolx| + |\bolx'|)\big)  
\big( |x_i-x_i'| + \frac{1}{\sqrt{N}}|\bolx - \bolx'| \big) \\
\leq & \; 4[H_b]_{_2}\big(1 + |\bolx| + |\bolx'|\big) |\bolx - \bolx'| \\
\big| G(\bar\mu(\bolx)) - G(\bar\mu(\bolx')) \big| & \leq \; \frac{[G]}{N} (|\bolx| + |\bolx'|)|\bolx - \bolx'|, 
\end{align} 
for all $x,x'$ $\in$ $\R^d$, $\bolx$ $=$ $(x_i)_{i\in\llbracket 1,N\rrbracket}$, $\bolx'$$=$ $(x'_i)_{i\in\llbracket 1,N\rrbracket}$ $\in$ $(\R^d)^N$, which yields by standard stability results for BSDE (see e.g. Theorems 4.2.1 and 5.2.1 in \cite{zha17}) that the function $v^N$ in \eqref{relYNv}  inherits the locally Lipschitz condition:  
\begin{align}
\big| v^N(t,\bolx) - v^N(t,\bolx') | & \leq \; C(1 + |\bolx| + |\bolx'|)|\bolx - \bolx'|, \quad \forall \bolx,\bolx' \in (\R^d)^N,  
\end{align}
for some constant $C$ (possibly depending on $N$). This implies 
\begin{align} \label{estiZ} 
|Z_t^{i,N}| & \leq \; 
C\big( 1 +  |\boX_t^N| \big),  \quad 0 \leq t \leq T, \quad i=1,\ldots,N, 
\end{align}
(this is clear when $v^N$ is smooth, and otherwise obtained by a mollifying argument as in Theorem 5.2.2 in \cite{zha17}).  
\end{Remark}

\begin{Remark}
Assumption \ref{assumption: solution} is verified in the case of linear quadratic control problems for which explicit smooth solutions are found in \cite{phawei18,phawei17} respectively without and with common noise. These papers prove that the second order Lions derivative $\partial^2_\mu$ is a continuous function of time which does not depend on the $\mu,x$ arguments hence is bounded whereas $\partial_\mu$ is affine in both the state and the first moment of the measure thus satisfies linear growth. Notice that in general, Assumptions \ref{assumption: Lipschitz} and \ref{assumption: linear} are not satisfied due to  the quadratic nature of $H_b$ in the $z$. 
However, in the uncontrolled case 
\begin{align}
    v(t,\mu) =\ & \E_{t,\mu}\Big[\int_t^T \Big( X_t^\top A(t) X_t + \E[X_t]^\top B(t) \E[X_t] +  C(t) X_t + D(t) \E[X_t] \Big) dt 
    \\ & \quad + X_T^\top E X_T + \E[X_T]^\top F \E[X_T] +  G X_T + H \E[X_T]\  \Big]\\
    \di X_t =\ & (b_0(t) + b_1(t)X_t + b_2(t) \E[X_t])\ \di t + \sigma(t)\ \di W_t, 
\end{align}
 we see that $v$ is a solution to the linear PDE
\begin{equation} 
\begin{cases}
 \partial_t v  +  \int_{\R^d} \big[ x^\top A(t) x +  \bar \mu^\top B(t)  \bar \mu +  C(t) x + D(t) \bar \mu & \\ + (b_0(t) + b_1(t)x + b_2(t) \bar \mu) \partial_\mu v(t,\mu) (x) & \\  + \frac{1}{2}{\rm tr}\big((\sigma\sigma\trans)(t)\partial_x\partial_\mu v(t,\mu)(x)\big) \big] \mu(\di x)  =     0, \quad \quad & (t,\mu) \in   [0,T)\times\Pc_2(\R^d),  \\
 v(T,\mu) = E\ \mathrm{Var}(\mu) + \bar \mu^\top (E+ F) \bar \mu  + (G+H) \bar \mu, \quad \quad & \mu \in \Pc_2(\R^d),
\end{cases}
\end{equation} where $\bar \mu = \int_{\R^d} x \mu(\di x)$, $\mathrm{Var}(\mu) = \int_{\R^d} x^2 \mu(\di x) - \bar \mu^2$. In that case both assumptions \ref{assumption: solution} and \ref{assumption: Lipschitz} are satisfied. 
\end{Remark}

\begin{Theorem} \label{th: N convergence}
Under Assumptions \ref{assumption: solution} and \ref{assumption: Lipschitz}, we have $\P$-almost surely
\begin{align}
\Ec_N^\mry  & \leq \;  \frac{C_y}{N}, 
\end{align} 
where $C_y$ $=$ $\frac{T}{2} e^{[H_b]_{1}T}  L \| \sigma \|^2_\infty$, with $\|\sigma\|_\infty$ $=$ $\sup_{(t,x,\mu) \in [0,T]\times\R^d\times\Pc_2(\R^d)} |\sigma(t,x,\mu)|$. 
\end{Theorem}

\begin{Theorem} \label{th: convergence Z}
Under Assumptions \ref{assumption: solution}, \ref{assumption: Lipschitz} and \ref{assumption: linear}, we have 
\begin{align}
\big\| \Ec_N^\mrz \big\|_{_2} & \leq \; \frac{C_z}{N^{\frac{1}{2}}}, 
\end{align} 
where $C_z$ $=$ $\|\sigma^+\|_\infty \sqrt{ 2T([H_1]_1 + [H_2]_1 L) \bar C_y^2 + \bar C_y T \|\sigma\|_\infty^2 L   + \frac{2}{c_0} \bar C_y^2 T \| H_2\|_\infty^2}$ and \\ $\bar C_y$ $=$ $\frac{T}{2} e^{([H_1]_1 + [H_2]_1 L)T}  L \| \sigma \|^2_\infty$. 
\end{Theorem}

\begin{Remark} \label{remerror}
Let us consider the global weak errors  on $v$ and its $L$-derivative $\partial_\mu v$ along the limiting McKean-Vlasov SDE, and defined by 
\begin{align}
E_N^\mry & := \;  \sup_{0\leq t\leq T} \big| \E[Y_t^N] - \E[ v(t,\P^0_{_{X_t}}) ] \big| \\
E_N^\mrz &:=  \frac{1}{N} \sum_{i=1}^N \Big( \int_0^T \Big| \E\big[NZ_t^{i,N}\big]  -  \E\big[ \partial_\mu v(t,\P^0_{_{X_t^i}}) (X_t^i) \big] \Big|^2 \di t \Big)^{\frac{1}{2}},
\end{align}
where $X^i$ has the same law than $X$, and with McKean-Vlasov dynamics as in \eqref{SDEMKV} but driven by $W^i$, $i$ $=$ $1,\ldots,N$. 
Then, they can be decomposed as 
\begin{align}
E_N^{\mry} & \leq \;  \E\big[\Ec_N^{\mry}\big] + \tilde E_N^{\mry}, \quad   E_N^{\mrz} \;  \leq \;  \big\| \Ec_N^{\mrz} \big\|_{_2}   + \tilde E_N^{\mrz},
\end{align}
where  $\tilde E_N^{\mry}$,  $\tilde E_N^{\mrz}$ are  the (weak)  propagation of chaos errors defined by 
\begin{align}
\tilde E_N^{\mry} & :=  \;  \sup_{0\leq t \leq T} \big| \E[ v(t,\bar\mu(\boX_t^N)) ] - \E[v(t,\P^0_{_{X_t}})] \big| \\
\tilde E_N^{\mrz} & := \; \frac{1}{N} \sum_{i=1}^N \Big(  \int_0^T \Big| \E\big[\partial_\mu v(t,\bar\mu(\boX_t^N))(X_t^{i,N})\big] - \E\big[ \partial_\mu v(t,\P^0_{_{X_t^i}})(X_t^i)\big] \Big|^2 \di t \Big)^{\frac{1}{2}}, 
\end{align} 
From the conditional propagation of chaos result, which states that for any fixed $k$ $\geq$ $1$, 
the law of $(X_t^{i,N})_{t\in [0,T]}^{i\in \llbracket 1,k\rrbracket}$ converges toward the conditional law of  $(X_t^{i})_{t\in [0,T]}^{i\in \llbracket 1,k\rrbracket}$, as $N$ goes to infinity, we deduce that 
$\tilde E_N^\mry,\tilde E_N^\mrz$ $\rightarrow$ $0$. Furthermore, under additional assumptions on $v$, we can obtain a rate of convergence. Namely, if $v(t,.)$ is Lipschitz uniformly in $t$ $\in$ $[0,T]$, with Lipschitz constant $[v]$, 
we have
\begin{align}
\tilde E_N^{\mry} & \leq \; [v] \sup_{0\leq t \leq T} \Big(\E[ \Wc_2(\bar\mu(\boX_t^N),\P^0_{_{X_t}})^2 ] \Big)^{\frac{1}{2}} \; =  \; O \Big(  N^{-\frac{1}{\max(d,4)}}\sqrt{ 1 + \ln(N) 1_{d=4}} \Big),   \nonumber \\
\mbox{ hence } \quad E_N^{\mry} & = \; O \Big(  N^{-\frac{1}{\max(d,4)}}\sqrt{ 1 + \ln(N) 1_{d=4} } \Big),  \label{rateEny}
\end{align} 
where we use the rate of convergence of empirical measures in Wasserstein distance stated in \cite{fougui15} (see also Theorem 2.12 in \cite{cardel2}),  
and since we have the standard estimate\\ $\E[\sup_{0\leq t\leq T}|X_t|^{4q}]$ $\leq$ $C(1+\|\mu_0\|^{4q}_{_{2q}})$ by Assumption \ref{assumption: Lipschitz}(iii).
The rate of convergence in \eqref{rateEny} is consistent with the one found in Theorem 6.17 \cite{cardel2} for  mean-field control problem. 
Furthermore, if the function $\partial_\mu v(t,.)(.)$ is Lipschitz in $(x,\mu)$ uniformly in $t$, then by the rate of convergence in Theorem 2.12 of  \cite{cardel2}, we have
\begin{align}
\tilde E_N^{\mrz} &  =  \; O \big( N^{-\frac{1}{\max(d,4)}}\big( 1 + \ln(N) 1_{d=4} \big)  \Big),   \; \mbox{ hence } \quad E_N^{\mrz}  =  \; O \big( N^{-\frac{1}{\max(d,4)}}\big( 1 + \ln(N) 1_{d=4} \big)  \Big). 
\end{align} 
\end{Remark} 

\begin{Remark}[Comparison with \cite{ganmayswi20}]\label{rem: comparison gangbo mayorga swiech}
In the related paper \cite{ganmayswi20} the authors consider a pure common noise case, that is $\sigma = 0 $ and restrict themselves to $\sigma_0(t,x_i,\bar\mu(\bolx))=\kappa I_d $ for $\kappa\in\R$.  If we consider these assumptions in our smooth setting, we directly see that $\Delta Y_s^N = 0 $ and $\Delta Z_s^N = 0 $ $ \P$ a.s.  Indeed by \eqref{BSDEMKVN} and \eqref{BSDEMKVNtilde} we notice that $(Y_t^N,Z_t^N)$ and \begin{equation}
    (\tilde Y_t^N  :=  v(t,\bar\mu(\boX_t^N)),  \{ \tilde Z_t^{i,N}  := \; \frac{1}{N} \partial_\mu v(t,\bar\mu(\boX_t^{N}))(X_t^{i,N}) , 
\;\;\;  i=1,\ldots,N\}),
\end{equation} solve the same BSDE therefore by existence and pathwise uniqueness for Lipschitz BSDEs the result follows. Moreover, \cite{ganmayswi20} does not allow $H$ to depend on $y$. 
Our approach allows to extend their findings to the case of idiosyncratic noises and in contrast to them we are able to choose a state-dependent volatility coefficient. Moreover we provide a convergence rate for the solution. However, we have to assume existence of a smooth solution for the master equation which is a restrictive assumption. 
\end{Remark}


\section{Proof of main results}  \label{secproof}

\subsection{Proof of Theorem \ref{th: N convergence}} 
 
\noindent  \underline{{\it  Step 1.}} 
Under the smoothness condition on $v$ in Assumption \ref{assumption: solution}, one can apply the standard It\^o's formula in $(\R^d)^N$ to the process $\tilde v^N(t,\boX_t^N)$ $=$ $v(t,\bar\mu(\boX_t^N))$, 
and get
\begin{align} \label{vNito} 
\tilde v^N(t,\boX_t^N) &= \;  \tilde v^N(T,\boX_T^N) - \int_t^T \partial_t  \tilde v^N(s,\boX_s^N) \di  s \\
& \quad  - \int_t^T B_N(s,\boX_s^N).D_\bolx \tilde v^N(s,\boX_s^N)  + \frac{1}{2}{\rm tr}\big(\Sigma_N(s,\boX_s^N) D_\bolx^2 \tilde v^N(s,\boX_s^N)\big)    \di s \\
& \quad - \;  \sum_{i=1}^N \int_t^T \big( D_{x_i}  \tilde v^N(s,\boX_s^{N}) \big)\trans  \sigma(s,X_s^{i,N}, \bar\mu(\boX_s^{N})) \di W_s^i \\
&  \quad - \;   \sum_{i=1}^N \int_t^T \big( D_{x_i}  \tilde v^N(s,\boX_s^{N}) \big)\trans    \sigma_0(s,X_s^{i,N}, \bar\mu(\boX_s^{N})) \di W_s^0,
\end{align}

Now, by setting (recall  \eqref{derivN}):  
\begin{align}
\tilde Y_t^N & := \; v(t,\bar\mu(\boX_t^N)) \; = \; \tilde v^N(t,\boX_t^N), \\
\tilde Z_t^{i,N} & := \; \frac{1}{N} \partial_\mu v(t,\bar\mu(\boX_t^{N}))(X_t^{i,N}) \; = \; D_{x_i} \tilde v^N(t,\boX_t^N), 
\;\;\;  i=1,\ldots,N,  \quad 0 \leq t\leq T,  
\end{align}
and using the relation \eqref{tildevN} satisfied by $\tilde v^N$ into 
\eqref{vNito}, we have for all $0\leq t\leq T$, 
 \begin{align} \label{BSDEMKVNtilde}
\tilde Y_t^N & = \; G\big(\bar\mu(\boX_{T}^N)\big) + \frac{1}{N} \sum_{i=1}^N \int_t^T H_b(s,X_s^{i,N},\bar\mu(\boX_{s}^N) ,\tilde Y_s^N, N \tilde  Z_s^{i,N}) \di s \\
& \quad \quad - \;  \frac{1}{2N^2} \sum_{i=1}^N \int_t^T {\rm tr}\big(\Sigma(s,X_s^{i,N}, \bar\mu(\boX_s^{N})) \partial^2_\mu v(s,\bar\mu(\boX_s^{N}))(X_s^{i,N},X_s^{i,N}) \big) \di s  \\
& \quad \quad  - \; \sum_{i=1}^N \int_t^T (\tilde Z_s^{i,N})\trans \sigma\big(s,X_s^{i,N},\bar\mu(\boX_{s}^N) \big) \di W_s^i 
- \; \sum_{i=1}^N \int_t^T (\tilde Z_s^{i,N})\trans \sigma_0\big(s,X_s^{i,N},\bar\mu(\boX_{s}^N) \big) \di W_s^0.  
\end{align}

\noindent  {\it \underline{Step 2}: Linearization}.  We set 
\begin{align}
\Delta Y_t^N &:= \; Y_t^N - \tilde Y_t^N, \quad \Delta Z_t^{i,N} :=  N(Z_t^{i,N}-\tilde Z_t^{i,N}), \;\;\;  i=1,\ldots,N,  \quad 0 \leq t\leq T, 
\end{align}
so that by \eqref{BSDEMKVN}-\eqref{BSDEMKVNtilde}, 
\begin{align} 
\Delta Y_t^N & = \;  \frac{1}{N} \sum_{i=1}^N \int_t^T \big[ H_b(s,X_s^{i,N},\bar\mu(\boX_{s}^N) ,Y_s^N, N Z_s^{i,N})  - H_b(s,X_s^{i,N},\bar\mu(\boX_{s}^N) ,\tilde Y_s^N, N \tilde Z_s^{i,N})\big]  \di s \\
& \quad \quad + \;  \frac{1}{2N^2} \sum_{i=1}^N \int_t^T {\rm tr}\big(\Sigma(s,X_s^{i,N}, \bar\mu(\boX_s^{N})) \partial^2_\mu v(s,\bar\mu(\boX_s^{N}))(X_s^{i,N},X_s^{i,N}) \big) \di s  \\
& \quad \quad  - \;  \frac{1}{N}  \sum_{i=1}^N \int_t^T (\Delta Z_s^{i,N})\trans \sigma\big(s,X_s^{i,N},\bar\mu(\boX_{s}^N) \big) \di W_s^i \\
& \quad \quad  - \;  \frac{1}{N}  \sum_{i=1}^N \int_t^T (\Delta Z_s^{i,N})\trans \sigma_0\big(s,X_s^{i,N},\bar\mu(\boX_{s}^N) \big) \di W_s^0,  \quad \quad   0 \leq t \leq T. \label{eq: error expression Y}
\end{align}
We now use the linearization method for BSDEs and rewrite the above equation as
\begin{align}
\Delta Y_t^N & =  \int_t^T \alpha_s \Delta Y_s^N \di s + \frac{1}{N} \sum_{i=1}^N \int_t^T  \beta_s^i.\Delta Z_s^{i,N}  \di s  \\
& \quad  -  \frac{1}{N}  \sum_{i=1}^N \int_t^T (\Delta Z_s^{i,N})\trans \sigma\big(s,X_s^{i,N},\bar\mu(\boX_{s}^N) \big) \di W_s^i   \\
& \quad  -  \frac{1}{N}  \sum_{i=1}^N \int_t^T (\Delta Z_s^{i,N})\trans \sigma_0\big(s,X_s^{i,N},\bar\mu(\boX_{s}^N) \big) \di W_s^0  \nonumber  \\
& \quad \; + \;  \frac{1}{2N^2} \sum_{i=1}^N \int_t^T {\rm tr}\big(\Sigma(s,X_s^{i,N}, \bar\mu(\boX_s^{N})) \partial^2_\mu v(s,\bar\mu(\boX_s^{N}))(X_s^{i,N},X_s^{i,N}) \big) \di s,  \label{linY} 
\end{align}
with 
\begin{align}
\begin{cases}
 &\alpha_s = \frac{1}{N} \Sum_{i=1}^N \frac{H_b(s,X_s^{i,N},\bar\mu(\boX_s^{N}), Y_s^N, N \tilde Z_s^{i,N}) - H_b(s,X_s^{i,N},\bar\mu(\boX_{s}^N),\tilde Y_s^N, N \tilde Z_s^{i,N})}{\Delta Y_s^N} \mathds{1}_{\Delta Y_s^N \neq 0} \\
 & \beta_s^i = \frac{H_b(s,X_s^{i,N},\bar\mu(\boX_s^{N}),Y_s^N, N Z_s^{i,N}) - H_b(s,X_s^{i,N},\bar\mu(\boX_s^{N}), Y_s^N, N \tilde Z_s^{i,N})}{|\Delta Z_s^{i,N}|^2} \Delta Z_s^{i,N} \mathds{1}_{\Delta Z_s^{i,N} \neq 0} 
 \end{cases} \label{eq: alpha}
\end{align}
for $i$ $=$ $1,\ldots,N$, and we notice by  Assumption \ref{assumption: Lipschitz}(iv) that the processes $\alpha$ and $\beta^i$ are bounded by $[H_b]_{_1}$.   Under Assumption \ref{assumption: Lipschitz}(ii), let us define the bounded processes 
$\lambda^i_s$ $=$ $\sigma^{+}(s,X_s^{i,N},\bar\mu(\boX_{s}^N))\beta_s^i$, $s$ $\in$ $[0,T]$, $i$ $=$ $1,\ldots,N$, and introduce the change of probability measure $\P^\lambda$ with Radon-Nikodym density: 
\begin{align}
\frac{\di \P^\lambda}{\di\P} &= \; \Ec_T^\lambda \; := \;  \exp\Big(\sum_{i=1}^N \int_0^T \lambda_s^i  \di W_s^i - \sum_{i=1}^N \frac{1}{2} \int_0^T |\lambda_s^i|^2 \di s\Big), 
\end{align}
so that by Girsanov's theorem: $\widetilde W_t^i$ $=$ $W_t^i - \int_0^t \lambda_s^i \di s$, $i$ $=$ $1,\ldots,N$,  and $W^0$ are independant Brownian motion under $\P^\lambda$. 
By applying It\^o's Lemma to $e^{\int_0^s  \alpha_s \di s}\Delta Y_t^N$ under $\P^\lambda$, we then obtain
\begin{align}
 \Delta Y_t^N &= \;  
 \frac{1}{2N^2} \sum_{i=1}^N \int_t^T e^{\int_t^s \alpha_u \di u} {\rm tr}\big(\Sigma(s,X_s^{i,N}, \bar\mu(\boX_s^{N})) \partial^2_\mu v(s,\bar\mu(\boX_s^{N}))(X_s^{i,N},X_s^{i,N}) \big) \di s \nonumber \\
 & \quad  -  \frac{1}{N}  \sum_{i=1}^N \int_t^T e^{\int_t^s \alpha_u \di u} (\Delta Z_s^{i,N})\trans \sigma\big(s,X_s^{i,N},\bar\mu(\boX_{s}^N) \big) \di \widetilde W_s^i \\
 & \quad  -  \frac{1}{N}  \sum_{i=1}^N \int_t^T e^{\int_t^s \alpha_u \di u} (\Delta Z_s^{i,N})\trans \sigma_0\big(s,X_s^{i,N},\bar\mu(\boX_{s}^N) \big) \di W_s^0, \quad 0 \leq t \leq T.  \label{eq: linearized Y}
\end{align}

\noindent {\it  \underline{Step 3}.} Let us check that the stochastic integrals in \eqref{eq: linearized Y}, namely $\int \tilde\Zc_s^{i,N}. \di \widetilde W_s^i$,  and 
$\int \tilde\Zc_s^{0,i,N}. \di W_s^0$ are  ``true"  martingales 
under $\P^\lambda$, where  
$\tilde Z_s^{i,N}$ $:=$  $e^{\int_0^s \alpha_u \di u}  \sigma\trans \big(s,X_s^{i,N},\bar\mu(\boX_{s}^N) \big) \Delta Z_s^{i,N}$, 
$\tilde Z_s^{0,i,N}$ $:=$  $e^{\int_0^s \alpha_u \di u}  \sigma_0\trans \big(s,X_s^{i,N},\bar\mu(\boX_{s}^N) \big) \Delta Z_s^{i,N}$,  $i$ $=$ $1,\ldots,N$, $0\leq s \leq T$. 

Indeed, for fixed $i$ $\in$ $\llbracket 1,N\rrbracket$, recalling that $\alpha$ is bounded, and  by the linear growth condition of  $\sigma_0$ from Assumption \ref{assumption: Lipschitz}(i),  we have
\begin{align}
\E^{\P^\lambda} \Big[ \int_0^T |\tilde Z_s^{0,i,N}|^2 \di s \Big] & \leq \; C  \E^{\P^\lambda} \Big[ \int_0^T  \big( |\sigma_0(s,0,\delta_0)|^2  + |X_s^{i,N}|^2 + \|\bar\mu(\boX_s^N)\|^2_{_2} \big) \\
& \quad \quad \quad \quad \quad  \quad  \big(  |NZ_s^{i,N}|^2 +  |\partial_\mu v(s,\bar\mu(\boX_s^{N}))(X_s^{i,N})|^2 \big) \di s \\
& \leq \; C \E\Big[ \Ec_T^\lambda \int_0^T N^2(|\sigma_0(s,0,\delta_0)|^4 + |\boX_s^N|^4)     \di s \Big] 
\end{align} 
where we use Bayes formula, the estimation \eqref{estiZ},  the growth condition on $\partial_\mu v(.)(.)$ in Assumption \ref{assumption: solution}, and noting that $|X_s^{i,N}|$ $\leq$ $|\boX_s^N|$, $\|\bar\mu(\boX_s^N)\|^2_{_2}$ $=$ $|\boX_s^N|^2/N$.  
By Hölder inequality with $q$ as in Assumption \ref{assumption: Lipschitz}(iii), and 
$\frac{1}{p} + \frac{1}{q} = 1$, the above inequality yields 
\begin{align}
\E^{\P^\lambda} \Big[ \int_0^T |\tilde Z_s^{0,i,N}|^2 \di s \Big] & \leq \; C N^2 \Big(\E\big[ |\Ec_T^\lambda|^p \big] \Big)^{\frac{1}{p}} \Big(  \E \Big[ \int_0^T (|\sigma_0(s,0,\delta_0)|^{4q}  + |\boX_s^N|^{4q} ) \di s \Big] \Big)^{\frac{1}{q}},
\end{align}
which is finite by \eqref{estiXiN}, and since $\lambda$ is bounded. This shows the square-integrable martingale property of $\int \tilde\Zc_s^{0,i,N}. \di W_s^0$ under $\P^\lambda$.  By the same arguments, we get 
the square-integrable martingale property of $\int \tilde\Zc_s^{i,N}. \di \tilde W_s^i$ under $\P^\lambda$

\vspace{1mm}

\noindent  {\it \underline{Step 4}: Estimation of $\Ec_N^\mry$.}
By taking the $\P^\lambda$ conditional expectation in \eqref{eq: linearized Y},  we obtain
\begin{align} \label{eq: expression of error Y under P lambda}
\Delta Y_t^N &=  
 \frac{1}{2N^2} \sum_{i=1}^N \E^{\P^\lambda} \Big[ \int_t^T e^{\int_t^s \alpha_u \di u} {\rm tr}\big(\Sigma(s,X_s^{i,N}, \bar\mu(\boX_s^{N})) \partial^2_\mu v(s,\bar\mu(\boX_s^{N}))(X_s^{i,N},X_s^{i,N}) \big) \di s \big| \Fc_t \Big],  
\end{align}
for all $t$ $\in$ $[0,T]$. Under the boundedness condition on $\Sigma$ $=$ $\sigma\sigma\trans$ in Assumption \ref{assumption: Lipschitz}(ii), and on $\partial_\mu^2 v$  in Assumption \ref{assumption: solution},  it follows immediately that 

\begin{align} \label{estimY} 
\Ec_N^\mry \; = \; \sup_{0\leq t \leq T} |\Delta Y_t^N| & \leq \;  \frac{C_y}{N}, 
 \quad a.s. 
\end{align}
where $C_y$ $=$ $\frac{T}{2} e^{[H_b]_{1}T}  L \| \sigma \|^2_\infty$, with $\|\sigma\|_\infty$ $=$ $\sup_{(t,x,\mu) \in [0,T]\times\R^d\times\Pc_2(\R^d)} |\sigma(t,x,\mu)|$.

\subsection{Proof of Theorem \ref{th: convergence Z}} \label{secproofZ} 

From \eqref{eq: alpha}, and under Assumption \ref{assumption: linear}(i) and (iii),  we see that 
\begin{align}
    \alpha_s & = \frac{1}{N} \Sum_{i=1}^N \frac{H_b(s,X_s^{i,N},\bar\mu(\boX_s^{N}), Y_s^N, N \tilde Z_s^{i,N}) - H_b(s,X_s^{i,N},\bar\mu(\boX_{s}^N),\tilde Y_s^N, N  \tilde Z_s^{i,N})}{\Delta Y_s^N} \mathds{1}_{\Delta Y_s^N \neq 0}\\
    & = \frac{1}{N} \Sum_{i=1}^N \frac{H_1(s,X_s^{i,N},\bar\mu(\boX_s^{N}), Y_s^N) - H_1(s,X_s^{i,N},\bar\mu(\boX_{s}^N),\tilde Y_s^N)}{\Delta Y_s^N} \mathds{1}_{\Delta Y_s^N \neq 0}\\
    & \;\;\; + \;  \frac{1}{N} \Sum_{i=1}^N  N \tilde Z_s^{i,N}. \frac{H_2(s,\bar\mu(\boX_s^{N}), Y_s^N) - H_2(s,\bar\mu(\boX_{s}^N),\tilde Y_s^N)}{\Delta Y_s^N} \mathds{1}_{\Delta Y_s^N \neq 0},\label{eq: alpha proof 2}
\end{align} is bounded by $[H_1]_1 + [H_2]_1 L$, recalling $N \tilde Z_s^{i,N}=\partial_\mu v(s,\bar\mu(\boX_{s}^N) )(X_s^{i,N})$. As a consequence, the proof of Theorem \ref{th: N convergence} still applies. Then by $\eqref{eq: error expression Y}$
\begin{align}
\Delta Y_t^N & =  \int_t^T \alpha_s \Delta Y_s^N \di s + \frac{1}{N} \sum_{i=1}^N \int_t^T  H_2(s,\bar\mu(\boX_s^{N}), Y_s^N). \Delta Z_s^{i,N}  \di s  \\
& \quad  -  \frac{1}{N}  \sum_{i=1}^N \int_t^T (\Delta Z_s^{i,N})\trans \sigma\big(s,X_s^{i,N},\bar\mu(\boX_{s}^N) \big) \di W_s^i   \\
& \quad  -  \frac{1}{N}  \sum_{i=1}^N \int_t^T (\Delta Z_s^{i,N})\trans \sigma_0(s,\bar\mu(\boX_s^{N})) \di W_s^0  \nonumber  \\
& \quad \; + \;  \frac{1}{2N^2} \sum_{i=1}^N \int_t^T {\rm tr}\big(\Sigma(s,X_s^{i,N}, \bar\mu(\boX_s^{N})) \partial^2_\mu v(s,\bar\mu(\boX_s^{N}))(X_s^{i,N},X_s^{i,N}) \big) \di s. 
\end{align}
By applying It\^o's formula to $|\Delta Y_t^N|^2$ in  \eqref{linY} under $\P$
\begin{align}
& |\Delta Y_0^N|^2 +  \frac{1}{N^2}  \int_0^T \sum_{i=1}^N \big|\sigma\trans\big(s,X_s^{i,N},\bar\mu(\boX_{s}^N) \big)  \Delta Z_s^{i,N} \big|^2 \di s 
+ |\sigma_0\trans(s,\bar\mu(\boX_s^{N}))\sum_{j=1}^N  \Delta Z_s^{j,N}|^2 \di s \\ & = 2 \int_0^T \alpha_s |\Delta Y_s^N |^2\di s + \frac{2}{N} \sum_{i=1}^N \int_0^T  \Delta Y_s^N H_2(s,\bar\mu(\boX_s^{N}), Y_s^N).\Delta Z_s^{i,N}  \di s  \\
& \quad  -  \frac{2}{N}  \sum_{i=1}^N \int_0^T \Delta Y_s^N(\Delta Z_s^{i,N})\trans \sigma\big(s,X_s^{i,N},\bar\mu(\boX_{s}^N) \big) \di W_s^i   \\
& \quad  -  \frac{2}{N}  \sum_{i=1}^N \int_0^T \Delta Y_s^N(\Delta Z_s^{i,N})\trans \sigma_0(s,\bar\mu(\boX_s^{N})) \di W_s^0  \nonumber  \\
& \quad \; + \;  \frac{1}{N^2} \sum_{i=1}^N \int_0^T \Delta Y_s^N{\rm tr}\big(\Sigma(s,X_s^{i,N}, \bar\mu(\boX_s^{N})) \partial^2_\mu v(s,\bar\mu(\boX_s^{N}))(X_s^{i,N},X_s^{i,N}) \big) \di s, 
\end{align}
so by taking expectation under $\P$, and  using the Cauchy-Schwarz inequality in $\R^d$
\begin{align}
&  \frac{1}{N^2}  \int_0^T \sum_{i=1}^N \E\Big[\big|\sigma\trans\big(s,X_s^{i,N},\bar\mu(\boX_{s}^N) \big)  \Delta Z_s^{i,N} \big|^2  
+ |\sigma_0\trans(s,\bar\mu(\boX_s^{N})) \sum_{j=1}^N  \Delta Z_s^{j,N}|^2 \Big] \di s\\ & \leq \E\Big[2 \int_0^T |\alpha_s| |\Delta Y_s^N |^2\di s + 2  \int_0^T  \Big|\Delta Y_s^N H_2(s,\bar\mu(\boX_s^{N}), Y_s^N).  \sum_{i=1}^N  \frac{\Delta Z_s^{i,N}}{N} \Big| \di s\Big]  \\
& \quad \; + \;  \frac{1}{N^2} \sum_{i=1}^N \int_0^T \E|\Delta Y_s^N{\rm tr}\big(\Sigma(s,X_s^{i,N}, \bar\mu(\boX_s^{N})) \partial^2_\mu v(s,\bar\mu(\boX_s^{N}))(X_s^{i,N},X_s^{i,N}) \big) |\di s,  \label{eq: err Z proof}
\\ & \leq \E\Big[2 \int_0^T |\alpha_s| |\Delta Y_s^N |^2\di s + 2  \int_0^T  \Big|\Delta Y_s^N H_2(s,\bar\mu(\boX_s^{N}), Y_s^N)\Big| \Big|\sum_{i=1}^N  \frac{\Delta Z_s^{i,N}}{N}\Big| \di s\Big]  \\
& \quad \; + \;  \frac{1}{N^2} \sum_{i=1}^N \int_0^T \E|\Delta Y_s^N{\rm tr}\big(\Sigma(s,X_s^{i,N}, \bar\mu(\boX_s^{N})) \partial^2_\mu v(s,\bar\mu(\boX_s^{N}))(X_s^{i,N},X_s^{i,N}) \big) |\di s,
\\ & \leq \E\Big[ \vartheta \int_0^T  |\Delta Y_s^N H_2(s,\bar\mu(\boX_s^{N}), Y_s^N)|^2 \di s + \frac{1}{\vartheta N^2} \int_0^T \Big|\sum_{i=1}^N \Delta Z_s^{i,N}\Big|^2  \di s\Big] + \frac{\tilde C_z}{N^2},
 \end{align} 
 by Young inequality for any $\vartheta >0$, boundedness of $\alpha$ (see \eqref{eq: alpha proof 2}), $\Sigma$, $\partial^2_\mu v$ and Theorem \ref{th: N convergence}, where 
 $\tilde C_z$ $=$ $2T([H_1]_1 + [H_2]_1 L) \bar C_y^2 + \bar C_y T \|\sigma\|_\infty^2 L$ and $\bar C_y$ $=$ $\frac{T}{2} e^{([H_1]_1 + [H_2]_1 L)T}  L \| \sigma \|^2_\infty$. Thus, by Assumption \ref{assumption: linear}(ii), and by choosing 
 $\vartheta = \frac{2}{c_0}$, it follows from the boundedness of $H_2$, and Theorem \ref{th: N convergence} that 
\begin{equation}
   \E\Big[\frac{1}{N} \sum_{i=1}^N \int_0^T|\sigma\trans\big(s,X_s^{i,N},\bar\mu(\boX_{s}^N) \big)\Delta Z_s^{i,N}|^2\ \di s \Big] 
   \; \leq \;  \frac{\tilde C_z + \frac{2}{c_0} \bar C_y^2 T \|H_2\|_\infty^2}{N}, 
\end{equation} 
which ends the proof by recalling that $\|\sigma^+ \|_\infty < + \infty$, using Cauchy-Schwarz inequality in $\R^N$ (in the form $\frac{1}{N}\sum_i^N \sqrt{|a_i|}\leq \sqrt{\frac{1}{N}\sum_i^N |a_i|}$) and Jensen inequality (in the form $\E[\sqrt{|X|}]\leq \sqrt{\E[|X|]}$).

\begin{small}

\printbibliography

\end{small}

\end{document}